\documentclass[12pt]{article}
\usepackage{graphicx}
\usepackage{amsmath}
\usepackage[dvips]{epsfig}
\usepackage{amssymb}

\makeatletter
\renewcommand{\section}{\@startsection
  {section}%
  {2}%
  {0mm}%
  {\baselineskip}%
  {0.5 \baselineskip}%
  {\centering}}
\makeatother
\begin{document}

\title {  A note on the generalized $q$-Bernoulli measures with weight $\alpha$ }
\author{D. V. Dolgy, T. Kim, S. H. Lee,  C. S.  Ryoo\\[0.5cm]
Institute of Mathematics and Computer Science, \\
Far Eastern National University, Vladivostok, 690060, Russia \\ \\
Division of General Education,  \\
 Kwangwoon University, Seoul 139-701,  Korea \\ \\
  Department of Mathematics, \\
          Hannam University, Daejeon 306-791, Korea \\[0.5cm]
 }

\date{}
\maketitle

 { \footnotesize{ \bf Abstract}\hspace{1mm}
In this paper we discuss new concept of the $q$-extension of
Bernoulli measure. From those measures, we derive some
interesting properties on the generalized $q$-Bernoulli numbers with weight $\alpha$ attached to $\chi$.}



\bigskip
\section{Introduction }
\bigskip

Let $p$ be a fixed   prime number.  Throughout this paper
$\mathbb{Z}_p$, $\Bbb Q_p$, and  $\mathbb{C}_p$ will,
respectively, denote the ring of $p$-adic rational integers, the
field of $p$-adic rational numbers, and the completion of
algebraic closure of $\Bbb Q_p$. Let $\Bbb N $ be the set of
natural numbers and $\Bbb Z_+ = \Bbb N \cup \{0 \}.$
 Let $\nu_p$
be the normalized exponential valuation of $\mathbb{C}_p$ with
$|p|_p=p^{-\nu_p(p)}=\dfrac{1}{p}$ (see [1-14]).

 When we talk of
$q$-extension, $q$ is variously considered as an indeterminate, a
complex number $q\in \mathbb{C},$ or a $p$-adic number $q\in\Bbb
C_p .$  Throughout this paper we assume that $q\in \mathbb{C}_p$
with  $|1-q|_p<1$ and we use the notation of  $q$-number as
$$[x]_q =\frac{1-q^x}{1-q}, \text{ (see [1-14])}.$$
Thus, we note that $\lim_{ q \rightarrow 1}[x]_q=x$.

In [2], Carlitz defined a set of numbers $\xi_k$ =$\xi_k(q)$
inductively by
$$\xi_{0}=1,   \quad  ( q  \xi+1)^k -  \xi_{k} =
\left \{\begin{array}{ll}
1, & \mbox{ if } k=1, \\
0, &  \mbox{ if } k>1,
\end{array} \right. \eqno(1)
$$
with the usual convention of replacing $\xi^{k}$ by $\xi_{k}.$

These numbers are $q$-extension of ordinary Bernoulli numbers
$B_k$. But they do not remain finite when $q=1$. So he modified
(1) as follows:

$$\beta_{0,q}=1,   \quad  q( q  \beta+1)^k -  \beta_{k, q} =
\left \{\begin{array}{ll}
1, & \mbox{ if } k=1, \\
0, &  \mbox{ if } k>1,
\end{array} \right. \eqno(2)
$$
with the usual convention of replacing $\beta^{k}$ by
$\beta_{k,q}.$

The numbers $\beta_{k,q}$ are called the $k$-th Carlitz
$q$-Bernoulli numbers.

In [1], Carlitz also considered  the extended Carlitz's
$q$-Bernoulli numbers as follows:

$$\beta_{0,q}^h=\dfrac{h}{[h]_q},   \quad  q^h( q  \beta^h+1)^k -  \beta_{k, q}^h =
\left \{\begin{array}{ll}
1, & \mbox{ if } k=1, \\
0, &  \mbox{ if } k>1,
\end{array} \right.
$$
with the usual convention of replacing $(\beta^h)^{k}$ by
$\beta_{k,q}^h.$

Recently, Kim considered $q$-Bernoulli numbers, which are
different extended Carlitz's $q$-Bernoulli numbers, as follows:
for $\alpha \in \Bbb N$ and $n \in \Bbb Z_+$,

$$\widetilde{\beta}_{0,q}^{(\alpha)}=1,   \quad  q( q^{\alpha
}  \widetilde{\beta}^{(\alpha)}+1)^n -  \widetilde{\beta}_{n,
q}^{(\alpha)} = \left \{\begin{array}{ll}
\dfrac{\alpha}{[\alpha]_q}, & \mbox{ if } n=1, \\ \\
0, &  \mbox{
if } n>1,
\end{array} \right. \eqno(3)
$$
with the usual convention of replacing
$(\widetilde{\beta}^{(\alpha)})^{k}$ by
$\widetilde{\beta}_{k,q}^{(\alpha)}$ (see [3]).

The numbers $\widetilde{\beta}_{k,q}^{(\alpha)}$ are called the
$k$-th $q$-Bernoulli numbers with weight $\alpha$.

For fixed $d \in \Bbb Z_+$ with  $(p,
d)=1$, we set
$$\aligned
&X= X_d=\varprojlim_N ( \mathbb{Z}/dp^N \mathbb{Z}), \quad X_1=\Bbb Z_p,\\
&X^*=\bigcup_{\begin{subarray}{l} 0<a<dp \\ (a,p)=1 \end{subarray}}( a+dp \mathbb{Z}_p),\\
&a+dp^N \mathbb{Z}_p=\{x\in X\mid x\equiv a\pmod{dp^N}\},
\endaligned$$
where $a\in \mathbb{ Z}$  satisfies the condition  $0\leq a<dp^N$.

Let $UD(\mathbb{Z}_p)$ be the space of uniformly differentiable
functions on $\mathbb{Z}_p.$
 For
$ f \in UD(\mathbb{Z}_p)$,
 the  $p$-adic $q$-integral on $\mathbb{Z}_p$ is  defined by Kim
 as follows:
$$I_q(f)=\int_{\mathbb{Z}_p}f(x) d\mu_{q}(x)=
\lim_{N \to \infty}\dfrac{1}{[p^N]_q}\sum_{x=0}^{ p^N-1}f(x)q^x,
\text{ (see [3, 4])}. \eqno(4)
$$
By  (3) and (4), the  Witt's formula for the $q$-Bernoulli numbers with weight $\alpha$
is given by
$$ \int_{\mathbb{Z}_p}[x]_{q^\alpha}^n  d\mu_{q}(x)
=\widetilde{\beta}_{n,q}^{(\alpha)}, \text{ where } n \in \Bbb
Z_+. \eqno(5)$$

The $q$-Bernoulli polynomials with weight $\alpha$ are also
defined by
$$ \widetilde{\beta}_{n,q}^{(\alpha)}(x)
= \sum_{l=0}^n \binom nl [x]_{q^{\alpha}}^{n-l} q^{ \alpha lx }
\widetilde{\beta}_{l,q}^{(\alpha)}. \eqno(6)$$

From (4), (5) and (6), we can derive the Witt's formula for
$\widetilde{\beta}_{n,q}^{(\alpha)}(x)$ as follows:
$$ \int_{\mathbb{Z}_p}[x+y]_{q^\alpha}^n  d\mu_{q}(y)
=\widetilde{\beta}_{n,q}^{(\alpha)}(x), \text{ where } n \in \Bbb
Z_+. \eqno(7)$$

For $n \in \Bbb Z_+ $ and $ d \in \Bbb N$, the distribution
relation  for the $q$-Bernoulli polynomials with weight $\alpha$
are known that
$$\widetilde{\beta}_{n,q}^{(\alpha)}(x)=\dfrac{[d]_{q^\alpha}^n}{[d]_q}
\sum_{a=0}^{d-1} q^a \widetilde{\beta}_{n,q^d}^{(\alpha)}\left(
\dfrac{x+a}{d} \right), \text{ (see [3])}. \eqno(8)$$
Recently, several authors have studied the $p$-adic $q$-Euler and Bernoulli measures on $\Bbb Z_p$ (see[8, 9, 11]).
The purpose of this paper is to construct $p$-adic $q$-Bernoulli
distribution with weight $\alpha$(= $p$-adic $q$-Bernoulli
unbounded measure with weight $\alpha$) on $\Bbb Z_p$ and to study their integral representations.
Finally, we construct the generalized $q$-Bernoulli numbers with weight $\alpha$ and investigate their properties
related to $p$-adic $q$-$L$-functions.

\bigskip
\section{  $p$-adic $q$-Bernoulli distribution with weight $\alpha$}
\bigskip

Let $X$ be any compact-open subset of $\Bbb Q_p$, such as $\Bbb
Z_p$ or $\Bbb Z_p^*$. A $p$-adic distribution $\mu$ on $X$ is
defined to be an additive map from the collection of compact open
set in $X$ to $ \Bbb Q_p$:

$$\mu \left(  \bigcup_{k=1}^n U_k \right) = \sum_{k=1}^n \mu (U_k)( \text{additivity}),$$
where $\{ U_1, U_2, \cdots, U_n \}$ is any collection of disjoint
compact opensets in $X$.

The set $\Bbb Z_p$ has a topological basis of compact open sets of
the form $ a+p^n \mathbb{Z}_p$.

Consequently, if $U$ is any compact open subset of $\Bbb Z_p$, it
can be written as a finite disjoint union of sets
$$U= \bigcup_{j=1}^k( a_j + p^N \mathbb{Z}_p), $$
where $N \in \Bbb Z_+$ and $a_1, a_2, \cdots, a_k \in \Bbb Z$ with
$0\leq a_i <p^N-1$

Indeed, the $p$-adic ball $ a+p^n \mathbb{Z}_p$ can be represented
as the union of smaller balls
$$a+p^n \mathbb{Z}_p =\bigcup_{b=0}^{p-1}( a+bp^n+p^{n+1} \Bbb Z_p).$$

\bigskip

{ \bf Lemma 1.} Every map $\mu$ from the collection of
compact-open sets in $X$ to $\Bbb Q_p$ for which
$$\mu(a+ p^N \Bbb Z_p) =\bigcup_{b=0}^{p-1}( a+bp^N+ d p^{N+1} \Bbb Z_p)$$
holds whenever $a+ p^N \Bbb Z_p \subset X$, extends to a $p$-adic
distribution on $X$.

\bigskip

Now we define a map $\mu_{k,q}^{(\alpha)}$ on the balls in $\Bbb
Z_p$ as follows:

$$\mu_{k,q}^{(\alpha)}(a+ p^n \Bbb Z_p)=\dfrac{[p^n]_{q^\alpha}^k}{[p^n]_q}
 q^a f_{k, q^{p^n}}^{(\alpha)} \left(
\dfrac{ \{ a \}_n }{p^n} \right), \eqno(9)$$ where $\{ a \}_n$ is
the unique number in the set $\{ 0,1, 2, \cdots, p^n-1\}$ such
that $\{ a \}_n \equiv a \pmod {p^n}$.

If $ a \in \{ 0,1, 2, \cdots, p^n-1\}$, then

$$  \aligned  & \sum_{b=0}^{p-1} \mu_{k,q}^{(\alpha)}(a+ b p^n+ p^{n+1} \Bbb Z_p)
 =  \sum_{b=0}^{p-1}\dfrac{[p^{n+1}]_{q^\alpha}^k}{[p^{n+1}]_q}
 q^{a+bp^n} f_{k, q^{p^{n+1}}}^{(\alpha)} \left(
\dfrac{  a+bp^n  }{p^{n+1}} \right)\\
  &= q^a \dfrac{[p^{n}]_{q^\alpha}^k}{[p^{n}]_q}\dfrac{[p]_{(q^{p^n})^\alpha}^k}{[p]_{q^{p^n}}} \sum_{b=0}^{p-1} q^{bp^n} f_{k, {(q^{p^n})^p}}^{(\alpha)}
  \left( \dfrac{  \frac{a}{p^n}+b  }{p} \right).
 \endaligned \eqno(10)
$$

From (10), we note that $\mu_{k,q}^{(\alpha)}$ is $p$-adic
distribution on $\Bbb Z_p$ if and only if
$$\dfrac{[p]_{(q^{p^n})^\alpha}^k}{[p]_{q^{p^n}}} \sum_{b=0}^{p-1} q^{bp^n} f_{k, {(q^{p^n})^p}}^{(\alpha)}
  \left( \dfrac{  \frac{a}{p^n}+b  }{p} \right)=f_{k, q^{p^{n}}}^{(\alpha)} \left(
\dfrac{  a }{p^{n}} \right). $$

\bigskip

{ \bf Theorem  2.} Let $\alpha \in \Bbb N$ and $k \in \Bbb Z_+$.
Then we see that $\mu_{k,q}^{(\alpha)}(a+ p^n \Bbb Z_p)$ is
$p$-adic distribution on $\Bbb Z_p$ if and only if
$$\dfrac{[p]_{(q^{p^n})^\alpha}^k}{[p]_{q^{p^n}}} \sum_{b=0}^{p-1}
q^{bp^n} f_{k, {(q^{p^n})^p}}^{(\alpha)}
  \left( \dfrac{  \frac{a}{p^n}+b  }{p} \right)=f_{k, q^{p^{n}}}^{(\alpha)} \left(
\dfrac{  a }{p^{n}} \right). $$

\bigskip

We set
$$ f_{k, q^{p^{n}}}^{(\alpha)} (x)=  \widetilde{\beta}_{k,q^{p^n}}^{(\alpha)}(x). \eqno(11)$$

From (9) and (11), we get
$$\mu_{k,q}^{(\alpha)}(a+ p^n \Bbb Z_p)= \dfrac{[p^n]_{q^\alpha}^k}{[p^n]_q}
 q^a  \widetilde{\beta}_{k,q^{p^n}}^{(\alpha)}\left(
\dfrac{  a }{p^{n}} \right). \eqno(12)$$

By (8), (12) and Theorem 2, we obtain the following theorem.

\bigskip

{ \bf Theorem  3.} Let $\mu_{k,q}^{(\alpha)}$ be given by
$$\mu_{k,q}^{(\alpha)}(a+ dp^N \Bbb Z_p)=\dfrac{[dp^N]_{q^\alpha}^k}{[dp^N]_q}
 q^a  \widetilde{\beta}_{k,q^{dp^N}}^{(\alpha)}\left(
\dfrac{  a }{ dp^{N}} \right). \eqno(13)$$
 Then  $\mu_{k,q}^{(\alpha)}$ extends to a $\Bbb Q(q)$--valued distribution on the compact open sets $U \subset X$.

\bigskip

From (13), we note that

$$  \aligned   \int_X  d\mu_{k,q}^{(\alpha)}(x) & = \lim_{N \rightarrow \infty}
\sum_{x=0}^{dp^N-1} \mu_{k,q}^{(\alpha)}(x+ d p^N \Bbb Z_p) \\
 & =  \lim_{N \rightarrow \infty}\dfrac{[dp^N]_{q^\alpha}^k}{[dp^N]_q}
 \sum_{a=0}^{dp^N-1}q^a\widetilde{\beta}_{k,q^{dp^N}}^{(\alpha)}\left(
\dfrac{  a }{ dp^{N}} \right) .
 \endaligned \eqno(14)
$$

By (8) and (14), we get

$$  \int_X  d\mu_{k,q}^{(\alpha)}(x) = \widetilde{\beta}_{k,q}^{(\alpha)} .
$$
Therefore, we obtain the following theorem.

\bigskip

{ \bf Theorem  4.} For $\alpha \in \Bbb N$ and $k \in \Bbb Z_+$,
we have  $$  \int_X  d\mu_{k,q}^{(\alpha)}(x) =
\widetilde{\beta}_{k,q}^{(\alpha)} .
$$
\bigskip

Let $\chi$ be Dirichlet character with conductor $d \in
\mathbb{N}$. Then we define the generalized $q$-Bernoulli numbers
attached to $\chi$ as follows:
$$  \aligned   \widetilde{\beta}_{n,\chi, q}^{(\alpha)} & = \int_X  \chi(x) [x]_{q^\alpha}^n d\mu_{q}(x)\\
&= \dfrac{[d]_{q^\alpha}^n}{[d]_q}
 \sum_{a=0}^{d-1} q^a \chi(a) \widetilde{\beta}_{n,q^{d}}^{(\alpha)}\left(
\dfrac{  a }{ d} \right).
\endaligned \eqno(15)$$

From (13) and (15), we can derive the following equation.

$$ \aligned  &  \int_X  \chi(x)  d\mu_{k,q}^{(\alpha)}(x) =  \lim_{N \rightarrow \infty}
\sum_{x=0}^{dp^N-1} \chi(x)  \mu_{k,q}^{(\alpha)}(x+ d p^N \Bbb Z_p) \\
 & =  \lim_{N \rightarrow \infty}\dfrac{[dp^N]_{q^\alpha}^k}{[dp^N]_q}
 \sum_{x=0}^{dp^N-1} \chi(x) q^x \widetilde{\beta}_{k,q^{dp^N}}^{(\alpha)}\left(
\dfrac{  x }{ dp^{N}} \right) \\
&= \dfrac{[d]_{q^\alpha}^k}{[d]_q}
 \sum_{a=0}^{d-1} q^a \chi(a) \left \{ \lim_{N \rightarrow \infty}\dfrac{[p^N]_{q^{\alpha d}}^k}{[p^N]_{q^d}}
 \sum_{x=0}^{p^N-1} q^{dx}  \widetilde{\beta}_{k,q^{dp^N}} \left(
\dfrac{  \dfrac{a}{d}+x }{ p^{N}} \right)  \right \}\\
&= \dfrac{[d]_{q^\alpha}^k}{[d]_q} \sum_{a=0}^{d-1} q^a \chi(a)
\widetilde{\beta}_{k,q^{d}}^{(\alpha)}\left( \dfrac{  a }{ d}
\right)= \widetilde{\beta}_{k, \chi, q}^{(\alpha)},
\endaligned$$
and

$$ \aligned  &  \int_{pX}  \chi(x)  d\mu_{k,q}^{(\alpha)}(x) =
  \lim_{N \rightarrow \infty}\dfrac{[dp^{N+1}]_{q^\alpha}^k}{[dp^{N+1}]_q}
 \sum_{x=0}^{dp^N-1} \chi(px) q^{px}  \widetilde{\beta}_{k,q^{dp^{N+1}}}^{(\alpha)}\left(
\dfrac{  px }{ dp^{N+1}} \right) \\ &=
\dfrac{[p]_{q^\alpha}^k}{[p]_q}
\dfrac{[d]_{q^{p\alpha}}^k}{[d]_{q^p}} \sum_{a=0}^{d-1} \chi(pa)
q^{pa}  \lim_{N \rightarrow \infty} \dfrac{[p^N]_{ q^{d p
\alpha}}^k}{ [p^N]_{ q^{d p }}} \sum_{x=0}^{p^N-1} q^{pdx}
\widetilde{\beta}_{k,q^{pdp^N}}^{(\alpha)}\left( \dfrac{p(xd+ a)
}{ pdp^N} \right)\\
& = \dfrac{[p]_{q^\alpha}^k}{[p]_q}
\dfrac{[d]_{q^{\alpha p}}^k}{[d]_{q^p}} \sum_{a=0}^{d-1} \chi(p)
\chi(a) q^{pa} \widetilde{\beta}_{k,
q^{pd}}^{(\alpha)}\left( \dfrac{  a }{ d} \right)
= \chi(p) \dfrac{[p]_{q^{\alpha}}^k}{[p]_{q}}
\widetilde{\beta}_{k, \chi, q^{p}}^{(\alpha)}.
\endaligned
$$

For $\beta ( \neq 1) \in X^*$, we have
$$ \int_{pX}  \chi(x)  d\mu_{k,q^{1/\beta}}^{(\alpha)}( \beta x)
= \chi (\dfrac{p}{\beta}) \dfrac{[p]_{q^{\alpha/\beta}}^k}{[p]_{
q^{ 1/\beta}}} \widetilde{\beta}_{k, \chi, q^{p/
\beta}}^{(\alpha)},
$$
and
$$ \int_{X}  \chi(x)  d\mu_{k,q^{1/\beta}}^{(\alpha)}( \beta x)
= \chi (\dfrac{1}{\beta})\widetilde{\beta}_{k, \chi, q^{1/
\beta}}^{(\alpha)}.
$$

Therefore, we obtain the following theorem.

\bigskip

{ \bf Theorem  5.} For $\beta ( \neq 1) \in X^*$, we have
$$ \aligned  &  \int_{X}  \chi(x)  d\mu_{k,q}^{(\alpha)}(x) = \widetilde{\beta}_{k, \chi, q}^{(\alpha)},\\
  &  \int_{pX}  \chi(x)  d\mu_{k,q}^{(\alpha)}(x) =\chi(p) \dfrac{[p]_{q^{\alpha}}^k}{[p]_{q}}
\widetilde{\beta}_{k, \chi, q^{p}}^{(\alpha)},\\
&  \int_{pX}  \chi(x)  d\mu_{k,q^{1/\beta}}^{(\alpha)}( \beta x) =
\chi (\dfrac{p}{\beta}) \dfrac{[p]_{q^{\alpha/\beta}}^k}{[p]_{ q^{
1/\beta}}} \widetilde{\beta}_{k, \chi, q^{p/ \beta}}^{(\alpha)},\\
&  \int_{X}  \chi(x)  d\mu_{k,q^{1/\beta}}^{(\alpha)}( \beta x) =
\chi (\dfrac{1}{\beta})   \widetilde{\beta}_{k, \chi, q^{1/
\beta}}^{(\alpha)}.
\endaligned
$$
Define
$$ \mu_{k, \beta, q}^{(\alpha)}(U)
=\mu_{k,  q}^{(\alpha)}(U)- \beta^{-1} \dfrac{[\beta^{-1}]_{
q^\alpha}^k}{[\beta^{-1}]_q}\mu_{k, q^{1/ \beta}}^{(\alpha)}(
\beta U). \eqno(16)$$ By simple calculation, we get

$$ \aligned   & \int_{X^*}  \chi(x)  d\mu_{k,\beta,
q}^{(\alpha)}(x)\\ & = \int_{X}  \chi(x)  d\mu_{k,
q}^{(\alpha)}(x)- \beta^{-1} \dfrac{[\beta^{-1}]_{
q^\alpha}^k}{[\beta^{-1}]_q} \int_{pX}
\chi(x) \mu_{k, q^{1/ \beta}}^{(\alpha)}(x)\\
 &= \widetilde{\beta}_{k, \chi, q}^{(\alpha)}-\chi(p) \dfrac{[p]_{q^{\alpha}}^k}{[p]_{q}}
\widetilde{\beta}_{k, \chi, q^{p}}^{(\alpha)},
\endaligned \eqno(17)
$$
and

$$ \aligned    & \dfrac{[\beta^{-1}]_{
q^\alpha}^k}{[\beta^{-1}]_q^k}  \int_{X^*}  \chi(x) d\mu_{k, q^{1/
\beta}}^{(\alpha)}(\beta x)
 =  \dfrac{[\frac{1}{\beta}]_{ q^\alpha}^k}{[\frac{1}{\beta}]_q}
\chi ({1}/{\beta})   \widetilde{\beta}_{k, \chi, q^{1/
\beta}}^{(\alpha)} - \chi ({p}/{\beta}) \dfrac{[
\frac{p}{\beta}]_{ q^\alpha}^k}{[\frac{p}{\beta}]_q}
\widetilde{\beta}_{k, \chi, q^{p/ \beta}}^{(\alpha)}.
\endaligned \eqno(18)
$$
By (16), (17) and (18), we get

$$ \aligned  & \int_{X^*}  \chi(x)  d\mu_{k,\beta, q}^{(\alpha)}(\beta
x)\\
& = \int_{X}  \chi(x)  d\mu_{k, q}^{(\alpha)}(x)- \beta^{-1}
\dfrac{[\beta^{-1}]_{ q^\alpha}^k}{[\beta^{-1}]_q} \int_{pX}
\chi(x) \mu_{k, q^{1/ \beta}}^{(\alpha)}( \beta x)\\
 &= \widetilde{\beta}_{k, \chi, q}^{(\alpha)}-\chi(p) \dfrac{[p]_{q^{\alpha}}^k}{[p]_{q}}
\widetilde{\beta}_{k, \chi, q^{p}}^{(\alpha)}-\dfrac{1}{\beta}
 \dfrac{[ \frac{1}{\beta}]_{
q^\alpha}^k}{[\frac{1}{\beta}]_q} \chi ({1}/{\beta})
\widetilde{\beta}_{k, \chi, q^{1/ \beta}}^{(\alpha)}\\
& \quad \quad  + \chi ({p}/{\beta}) \dfrac{[ \frac{p}{\beta}]_{
q^\alpha}^k}{[\frac{p}{\beta}]_q} \widetilde{\beta}_{k, \chi,
q^{p/ \beta}}^{(\alpha)}.
\endaligned \eqno(19)
$$

Now we define the operator $\chi^y=\chi^{y, k, \alpha :q}$ on
$f(q)$ by
$$ \chi^y f(q)=\chi^{y, k, \alpha :q} f(q)= \dfrac{[ y]_{
q^\alpha}^k}{[y]_q} \chi(y) f(q^y). \eqno(20)$$

Thus, by (20), we get

$$ \aligned  & \chi^{x, k, \alpha :q} \circ  \chi^{y, k, \alpha
:q}f(q)  =  \chi^{x, k, \alpha :q} \dfrac{[ y]_{
q^\alpha}^k}{[y]_q} \chi(y) f(q^y) \\
&=\dfrac{[ y]_{ q^\alpha}^k}{[y]_q} \chi(y)  \chi(x) \dfrac{[ y]_{
q^{\alpha y}}^k}{[y]_{q^y}} \chi(y) f(q^{xy})\\
&=\dfrac{[ xy]_{ q^\alpha}^k}{[xy]_q} \chi(xy)    f(q^{xy})\\
&= \chi^{xy, k, \alpha :q}f(q)= \chi^{xy} f(q).
\endaligned \eqno(21)
$$
Let us define $\chi^{x}  \chi^{y}=\chi^{x, k, \alpha :q} \circ
\chi^{y, k, \alpha :q}$. Then we have $$\chi^{x}
\chi^{y}=\chi^{xy}.$$ From the definition of $\chi^{x}$, we can
easily derive the following equation.
$$(1-\chi^p)\left( 1-\frac{1}{\beta} x^{1/\beta}\right)= 1-\frac{1}{\beta}x^{1/\beta}-\chi^p + \frac{1}{\beta}x^{p/\beta}. $$

Let $f(q)=\widetilde{\beta}_{k, \chi, q}^{(\alpha)}$. Then we get
$$ \aligned  &(1-\chi^p)\left(1-\frac{1}{\beta} x^{1/\beta}\right) \widetilde{\beta}_{k, \chi, q}^{(\alpha)}  \\
&=\widetilde{\beta}_{k, \chi, q}^{(\alpha)}-\frac{1}{\beta}
\dfrac{[ \frac{1}{\beta}]_{ q^\alpha}^k}{[\frac{1}{\beta}]_q}\chi
({1}/{\beta}) \widetilde{\beta}_{k, \chi, q}^{(\alpha)}- \dfrac{[
p]_{ q^\alpha}^k}{[p]_q} \chi (p) \widetilde{\beta}_{k, \chi,
q^{p}}^{(\alpha)}
 \\ & \quad \quad + \dfrac{1}{\beta} \dfrac{[
\frac{p}{\beta}]_{ q^\alpha}^k}{[\frac{p}{\beta}]_q}
 \chi ({p}/{\beta})  \widetilde{\beta}_{k, \chi, q^{p/ \beta}}^{(\alpha)}.
\endaligned \eqno(22)
$$

By (19) and (22), we obtain the following equation:
$$  \int_{X^*}  \chi(x)  d\mu_{k,\beta, q}^{(\alpha)}(\beta
x)=(1-\chi^p)\left(1-\frac{1}{\beta} x^{1/\beta}\right)
\widetilde{\beta}_{k, \chi, q}^{(\alpha)}, $$ where $\beta( \neq
1) \in X^*$.

 \end{document}